\documentclass[11pt,leqno,english]{amsart}
\makeatletter
\def\paragraph{\@startsection{paragraph}{4}%
\z@{.5\linespacing}{.5\linespacing}%
{\bfseries}}
\def\@secnumfont{\ifnum\@toclevel>1\bfseries\else\mdseries\fi}
\makeatother

\usepackage[left=1.3in,top=1.3in,right=1.3in,nofoot]{geometry}
\usepackage{amsmath}
\usepackage{amsfonts}
\usepackage{amssymb}
\usepackage[bookmarksopen,bookmarksdepth=2]{hyperref}

\newcommand{\Q}{\mathbb{Q}}

\DeclareMathOperator{\GL}{GL}
\DeclareMathOperator{\PGL}{PGL}
\DeclareMathOperator{\SL}{SL}

\DeclareMathOperator{\Ind}{{Ind}}
\DeclareMathOperator{\ind}{{ind}}

\DeclareMathOperator{\car}{char}

\newtheorem{question}{Question}

\usepackage{xcolor}
\definecolor{teal}{rgb}{0.0, 0.5, 0.5}
\definecolor{forest}{rgb}{0.13, 0.55, 0.13}

\begin{document}

\title{Questions on mod $p$ representations of reductive $p$-adic groups}

\author{N. Abe}
\address[N. Abe]{Department of Mathematics, Hokkaido University, Kita 10, Nishi 8, Kita-Ku, Sapporo, Hokkaido, 060-0810, Japan}
\thanks{The first-named author was supported by JSPS KAKENHI Grant Number 26707001.}
\email{abenori@math.sci.hokudai.ac.jp}

\author{G. Henniart} 
\address[G. Henniart]{Universit\'e de Paris-Sud, Laboratoire de Math\'ematiques d'Orsay, Orsay cedex F-91405 France;
CNRS, Orsay cedex F-91405 France}
\email{Guy.Henniart@math.u-psud.fr}

\author{F. Herzig} 
\address[F. Herzig]{Department of Mathematics, University of Toronto,
  40 St.\ George Street, Room 6290, Toronto, ON M5S 2E4, Canada}
\thanks{The third-named author was partially supported by a Sloan Fellowship and an NSERC grant.}
\email{herzig@math.toronto.edu}

\author{M.-F. Vign\'eras}
\address[M.-F. Vign\'eras]{Institut de Math\'ematiques de Jussieu, 175 rue du Chevaleret, Paris 75013 France}
\email{vigneras@math.jussieu.fr}
\maketitle

\bigskip 0) {\bf Introduction}

We compiled these questions for the workshop \emph{Geometric methods in the mod p local Langlands correspondence}
held in June 2016 at the Centro di Ricerca Matematica Ennio de Giorgi in Pisa. We thank the organizers, Michael Harris and Peter
Schneider, for inviting us.

\bigskip The following is a preliminary and a bit informal discussion of questions raised by our work on modulo $p$ admissible smooth representations of reductive $p$-adic groups, for which some kind of answer is required if we want to venture into derived algebraic geometry for further study.

\bigskip 1) {\bf Admissibility questions}

The framework is the following: $p$ is a prime number, $F$ is a finite extension of $\mathbb Q_p$ or of $\mathbb F_p((T))$, $G$ is a connected reductive $F$-group, $R$ is the coefficient field, which unless otherwise stated, is algebraically closed of characteristic $p$. We examine smooth $R$-representations of $G(F)$, where
a representation of $G(F)$ on an $R$-vector space $V$ is
\emph{smooth} if the $G(F)$-stabilizer of any vector $v\in V$ is open. Such a representation is
\emph{admissible} if moreover the subspace $V^J$ of fixed vectors fixed under any open subgroup $J$ of $G(F)$ has finite dimension; actually it is enough to require that for one open pro-$p$ subgroup of $G(F)$.

Our joint work \cite{ahhv} gives a complete classification of irreducible admissible $R$-representations of $G(F)$ in terms of supercuspidal $R$-representations of Levi subgroups of $G(F)$ -- where a supercuspidal representation is
an irreducible admissible representation which is not a subquotient of a representation obtained by parabolic induction of an irreducible admissible representation of a proper Levi subgroup.

The requirement of admissibility, both in the definition of supercuspidal and in the classification, is a bit awkward. Indeed, if $C$ is an algebraically closed field of characteristic different from $p$, it is known that an irreducible smooth $C$-representation of $G(F)$ is admissible.

\begin{question}\label{qu:1}
  Is any irreducible smooth $R$-representation of $G(F)$ admissible?
\end{question}
The answer is yes when $G(F)=\GL(2,\Q_p)$ \cite{MR3076828}. It is also yes when $G$ is anisotropic modulo its centre: in that case, $G(F)$ divided by its centre is compact and all smooth irreducible representations of $G(F)$ are finite-dimensional.

An affirmative answer to Question~\ref{qu:1} has a number of desirable consequences. In an
irreducible admissible $R$-representation $\pi$ of $G(F)$ the centre of $G(F)$ acts via a character, called the central character of $\pi$. The following is weaker than Question~\ref{qu:1}.

\begin{question}\label{qu:2}
  Does any irreducible smooth $R$-representation of $G(F)$ possess a central character?
\end{question}

The answer is yes if $R$ is uncountable but it is unknown if $R$ is an algebraic closure of $\mathbb F_p$, for example. Answering Question~\ref{qu:2} might be the first step in answering Question~\ref{qu:1}.

\begin{question}\label{qu:3}
  Does $G(F)$ possess supercuspidal $R$-representations?
\end{question}

All supercuspidal representations of $\GL(2,\mathbb Q_p)$ are known (\cite{BP}, building on \cite{bib:BL-general}). Otherwise, many are constructed, though not in an explicit way, when $G=\GL_2$ and $F/\mathbb Q_p$ is an unramified extension (\cite{MR2128381}, \cite{BP}) and in a few other low rank cases. A local-global construction yields such supercuspidal representations for $G = \GL_n$, provided $F$ has characteristic $0$.

\bigskip For the following questions, the answer may depend on the characteristic of $F$, $0$ or $p$.

When $F$ has characteristic $0$, it is known that a quotient of an admissible $R$-representation of $G(F)$ is still admissible. However, when $F$ has characteristic $p$, we can construct an admissible $R$-representation of $\mathbb G_m(F)=F^\times$ with a quotient which is not admissible.
That quotient has infinite length so we might ask:

\begin{question}\label{qu:4}
  Assume $\car F=p$. Let $V$ be a finite length admissible $R$-representation of $G(F)$. Is every quotient of $V$ admissible?
\end{question}
Again, this is weaker than Question~\ref{qu:1}.

\bigskip If $P$ is a parabolic $F$-subgroup of $G$ and $M$ a Levi component of $P$, the parabolic induction functor $\Ind_P^G$ from smooth $R$-representations of $M(F)$ to
smooth $R$-representations of $G(F)$ preserves admissibility. The functor has a left adjoint, the usual Jacquet functor taking coinvariants under the unipotent radical $N(F)$ of $P(F)$, and also a right adjoint. We can prove that the right adjoint respects admissibility. On the other hand, the Jacquet functor respects admissibility, provided $F$ has characteristic $0$ (\cite{bib:emerton-ordinary2} when $R$ is finite).
We did not yet check if this remains true when $R$ is an algebraically closed field of characteristic $p$.

\begin{question}\label{qu:5}
  Assume $\car F=p$. Let $V$ be a finite length admissible $R$-representation of $G(F)$. Is
    the representation $V_{N(F)}$ of $M(F)$ admissible?
\end{question}
In a different direction, let $\varphi:G'\to G$ be a central isogeny of connected reductive $F$-groups and let $Z_G$ the center of $G$. If $F$ has characteristic $0$, the quotient $G(F)/Z_G(F) G'(F)$ is finite, whereas it is only compact when $F$ has characteristic $p$.

\begin{question}\label{qu:6}
  Assume $\car F=p$. Let $V$ be a finite length admissible $R$-representation of $G(F)$. If
    we view $V$ as a representation of $G'(F)$ via $\varphi$, do we get a finite length $R$-representation of $G'(F)$ with
    admissible subquotients?
\end{question}
The first case to look at is that of $\SL(2)\to \PGL(2)$ when $p=2$.

When $\car F=0$, Jan Kohlhaase \cite{kohlhaase} has investigated contragredients for irreducible admissible representations of $G(F)$. In particular, he proved that the contragredient of such a representation $V$ is $0$ unless $V$ has finite dimension. We have extended that result to the case where $\car F=p$. But Kohlhaase went further to study the derived functors of the contragredient functor.

\begin{question}\label{qu:7}
 Assume $\car F=p$. Do the derived functors of the contragredient functor lead to some kind
    of duality for admissible representations of $G(F)$?
\end{question}

{\bf Aside questions} Before turning to questions centering on weights and eigenvalues, we mention some of our current explorations. If $R$ has characteristic $p$ but is not necessarily algebraically closed, can we still get a classification of irreducible admissible $R$-representations of $G(F)$? Is any such representation actually defined over a field of finite type over $\mathbb F_p$? Is any supercuspidal $R$-representation of $G(F)$, whose central character has finite order, definable over a finite extension of $\mathbb F_p$?

\bigskip 2) {\bf Weights and eigenvalues}

Our classification of irreducible admissible representations of $G(F)$ uses weights and eigenvalues, to which we now turn as they raise their own set of questions.

We choose a special parahoric subgroup $K$ of $G(F)$. Let
$V$ be an irreducible $R$-representation of $K$, $\ind_K^{G(F)}V$ the compactly induced (smooth) representation of $G(F)$, and let $\pi$ be an irreducible smooth $R$-representation of $G(F)$.

We say that $V$ is a \emph{weight} of $\pi$ if $\pi$ is isomorphic to a quotient of $\ind_K^{G(F)} V$.
Every irreducible smooth $R$-representation of $G(F)$ admits a weight.

Let $Z(G,K,V)$ be the center of the $R$-algebra $H(G,K,V)$ of $G(F)$-intertwiners of $\ind_K^{G(F)} V$. A homomorphism $\chi:Z(G,K,V)\to R$
is called a \emph{Hecke eigenvalue} of $V$ in $\pi$ if $\pi$ is isomorphic to a quotient of
$$\pi(V,\chi)= \chi \otimes_{Z(G,K,V)}\ind_K^{G(F)} V .$$
Every weight $ V$ of an irreducible admissible $R$-representation $\pi$ of $G(F)$ admits a Hecke eigenvalue.

\begin{question}\label{qu:8}
  Let $\pi$ be an irreducible smooth $R$-representation of $G(F)$. Is it true that $\pi$
    admits a Hecke eigenvalue $\chi$ for some weight $ V$?
\end{question}
If Question~\ref{qu:1} has a negative answer, Question~\ref{qu:8} could still have a positive one. Irreducible quotients of $\pi(V,\chi)$ are still amenable to the techniques of our work. On the other hand, one could try and answer Question~\ref{qu:1} positively by dealing first with Question~\ref{qu:8}, and then proving that the irreducible quotients of $\pi(V,\chi)$ are admissible. Note that $\pi(V,\chi)$ has a central character, and hence also its subquotients.

\bigskip In previous work, for a parabolic subgroup $P=MN$ of $G$ in good position with respect to $K$, we constructed an algebra homomorphism (called \emph{Satake homomorphism}) $Z(G,K,V)\to Z(M, K\cap M(F), V_{K\cap N(F)})$ -- here $V_{K\cap N(F)}$ is an irreducible representation of the parahoric subgroup $K\cap M(F)$ of $M(F)$. We defined a \emph{supersingular} $R$-representation of $G(F)$ to be an irreducible admissible representation $\pi$ such that for all weights $V$ of $\pi$ all Hecke eigenvalues $\chi : Z(G,K,V)\to R$ of $V$ in $\pi$ are supersingular, i.e.\ never factor through the Satake homomorphism when $P$ is a proper parabolic subgroup of $G$. We proved that if $\pi$ is an irreducible admissible representation, then $\pi$ is supersingular if and only if it is supercuspidal, and that if $\pi$ admits a supersingular Hecke eigenvalue, it is supersingular.

\begin{question}\label{qu:9}
  Let $\pi$ be an irreducible smooth $R$-representation of $G(F)$ admitting a supersingular
    Hecke eigenvalue. Is it true that all its Hecke eigenvalues are supersingular?
\end{question}

The kernel of the natural map $$\ind_K^{G(F)}V\to \prod_\chi \pi(V,\chi)$$ is a $Z(G,K,V)$-representation $I_1(G,K,V)$ of $G(F)$. By induction, for any $n\geq 1$, we define
$ I_{n+1}(G,K,V) $ as the kernel of the natural map $ I_{n}(G,K,V)\to \prod_\chi \pi_n(V,\chi)$ where $\pi_n(V,\chi)=\chi \otimes I_n(G,K,V)$.

\begin{question}\label{qu:10}
  Is the decreasing filtration $( I_n(G,K,V))$ of $\ind_K^{G(F)}V $ finite?
\end{question}

Let $I$ be a pro-$p$ Iwahori subgroup of $G(F)$, chosen in $K$ and in good position. Let $H$ be its Hecke algebra, that is the algebra of endomorphisms of the $R$-representation $\ind_{ I}^{G(F)}1_R $ compactly induced by the trivial representation of $I$ on $R$.
If $V$ is an irreducible smooth $R$-representation of $K$, the action of $H(G,K,V)$ on $\ind_K^{G(F)}V$ commutes with the action of $G(F)$. It is known (Ollivier, Vign\'eras) that the $H$-module of $I$-invariants is a cyclic $H$-module and that the action of $H(G,K,V)$ on $(\ind_K^{G(F)}V)^I$ induces an isomorphism of $H(G,K,V)$ onto the $H$-endomorphisms of $(\ind_K^{G(F)}V)^I$. The action of $H(G,K,V)$ and of the center of $H$ on $(\ind_K^{G(F)}V)^I$ are almost the same. As $H$ is a finite module over its center, the space $ \chi\otimes (\ind_K^{G(F)}V)^I$ is finite-dimensional for any Hecke eigenvalue $\chi$ of $V$.

\begin{question}\label{qu:11}
  Let $\chi$ be a Hecke eigenvalue for the weight $V$. Is the natural map $ \chi\otimes
    (\ind_K^{G(F)}V)^I \to \pi(V,\chi)$ injective?
\end{question}
The map $\chi\otimes (\ind_K^{G(F)}V)^I \to \pi(V,\chi)^I$ cannot be surjective when $ \pi(V,\chi)^I$ is infinite-dimensional. This is the case when $F$ is a proper unramified extension of $\mathbb Q_p$, $G = \GL(2 )$ and $\chi$ is supersingular (Breuil, Morra).

For any connected reductive $F$-group $G$, any special parahoric subgroup $K$ and any weight $V$, there exist supersingular
Hecke eigenvalues $\chi : Z(G,K,V)\to R$.

\begin{question}\label{qu:12}
  If $\chi : Z(G,K,V)\to R$ is supersingular, is $\pi(V,\chi)$ non-zero?
\end{question}

Note that Question~\ref{qu:1} and Question~\ref{qu:12} (for one choice of $(V,\chi)$) together imply Question~\ref{qu:3}.  To
see this, pick $(V,\chi)$ such that $\chi$ is supersingular and $\pi(V,\chi)$ is non-zero. The representation $\pi(V,\chi)$
admits an irreducible quotient $\pi$, since it is finitely generated.  If Question~\ref{qu:1} holds, or in fact even if
Question~\ref{qu:1} holds only for irreducible representations with a central character, then $\pi$ is admissible, and $\pi$
is supercuspidal by our results discussed above.

(In \cite[\S V]{ahhv} we defined a non-zero quotient $R_G$ of $Z(G,K,V)$ which has the property that a character
$\chi : Z(G,K,V)\to R$ is supersingular if and only if $\chi$ factors through $R_G$, and we showed that
$\pi(V) := R_G \otimes_{Z(G,K,V)}\ind_K^{G(F)} V$ is a free $R_G$-module. It follows that Question~\ref{qu:12}
is equivalent to asking whether $\pi(V) \ne 0$. In particular, the answer to Question~\ref{qu:12} is independent of
$\chi$, for any fixed weight $V$.)

The following question may be related. There is a definition of what it means for a finite-dimensional $H$-module to be supersingular. 

\begin{question}
  If $X$ is a supersingular $H$-module, is $X\otimes_H \ind_{I}^{G(F)} 1_R$ non-zero?
\end{question}

\bigskip {\bf Aside question} If $\pi$ is an irreducible admissible $R$-representation of $G(F)$ it is equivalent to say that
$\pi$ is supersingular or that the space of $I$-fixed vectors is supersingular as $H$-module (Ollivier, Vign\'eras). Can we extend this
equivalence to the case where $\pi$ is only an irreducible quotient of some $\pi(V,\chi)$?

\bibliography{Pise}
\bibliographystyle{amsalpha}

\end{document}